\documentclass[12pt]{article}
\usepackage{amssymb,latexsym, amsmath, enumerate}
\usepackage{color}
\usepackage{hyperref}
\hypersetup{
                CJKbookmarks=true
	colorlinks=true,
	linkcolor=blue,
	filecolor=blue,      
	urlcolor=blue,
	citecolor=blue,
}

\def\R{\mathbb R}
\def\N{\mathbb N}
\def\C{\mathbb C}
\def\H{\mathbb H}

\newtheorem{theorem}{Theorem}
\newtheorem{lemma}{Lemma}
\newtheorem{remark}{Remark}

\newtheorem{corollary}{Corollary}
\newtheorem{definition}{Definition}
\newtheorem{example}{Example}

\title{On the method of differentiation and its application to asymptotics for the heat kernel on H-type groups}
\author{Ye Zhang}
\date{}
\begin{document}

\renewcommand{\theequation}{\thesection.\arabic{equation}}
\setcounter{equation}{0} \maketitle

\vspace{-1.0cm}

\bigskip

{\bf Abstract.} The aim of this note is twofold. The first one is to find conditions on the asymptotic sequence which ensures
differentiation of a general asymptotic expansion with respect to it. %such an asymptotic sequence.
Our method results from the classical one but generalizes it. As an application, our second aim is to give sharp asymptotic
estimates at infinity of the heat kernel on H-type groups by the method of
differentiation provided we have the result of the isotropic Heisenberg groups.

\medskip


\medskip


\medskip

\section{Introduction}
\setcounter{equation}{0}
It is well-known that if $f$ is analytic in an open sector $\Delta$ in the complex plane $\C$ (for example, $\Delta = \{z: |z| > R > 0, |\arg(z)|< \theta_0 < \pi\}$) and
\begin{align*}
f(z) \sim \sum_{n=0}^{\infty} a_n z^{-n}, \qquad \mbox{$z \to \infty$ in $\Delta$,}
\end{align*}
that is, for any $N \ge 0$, we have
\begin{align*}
f(z) - \sum_{n=0}^N a_n z^{-n} = o(z^{-N}), \qquad \mbox{$z \to \infty$ in $\Delta$,}
\end{align*}
then in any closed sector $\Delta^{\prime}$ of $\Delta$ we have
\begin{align*}
f^{\prime}(z) \sim - \sum_{n=1}^{\infty} n a_n z^{-n-1}, \qquad \mbox{$z \to \infty$ in $\Delta^{\prime}$.}
\end{align*}
For this classical result, one can refer to \cite[p. 21]{O74} or \cite[Theorem 4, p. 8]{W89}.
One aim of this note is to obtain a similar result for a general asymptotic expansion with respect to a particular kind of asymptotic sequence.
More precisely, we give some conditions on the asymptotic sequence $\{g_n\}$ as $z \to \infty$ in $\Delta$ (that is, a sequence of functions such that for any $n \ge 0$ we have $g_{n+1} = o(g_n)$ as $z \to \infty$ in $\Delta$) such that if $f$ and ${f_n}$ are analytic in $\Delta$ and
\begin{align*}
f(z) \sim \sum_{n=0}^{\infty} f_n(z), \qquad \{g_n\},\qquad \mbox{$z \to \infty$ in $\Delta$,}
\end{align*}
that is, for any $N \ge 0$ we have
\begin{align*}
f(z) -  \sum_{n=0}^N f_n(z) = o(g_N(z)), \qquad \mbox{$z \to \infty$ in $\Delta$,}
\end{align*}
then in some closed sector $\Delta^{\prime}$ of $\Delta$ we have $\{g_n^{\prime}\}$ is an asymptotic sequence as $z \to \infty$ in $\Delta^{\prime}$ and
\begin{align*}
f^{\prime}(z) \sim \sum_{n=0}^{\infty} f_n^{\prime}(z), \qquad \{g_n^{\prime}\},\qquad \mbox{$z \to \infty$ in $\Delta^{\prime}$.}
\end{align*}
For the rest of this note we show an application of the method of differentiation to the asymptotic behaviour at infinity of the heat kernel on H-type groups provided we have the corresponding result of the isotropic Heisenberg groups.
Asymptotic estimates at infinity of the heat kernel on H-type groups have been studied well in the last decades.
For example, in \cite{E09}, Eldredge
provides precise upper and lower bounds of the heat kernel on H-type groups and in \cite{Li10}, Li provides asymptotic estimates at infinity of the heat kernel
on H-type groups.

As a result, gradient bound estimates of the heat semigroup, that is,
\begin{align} \label{HGE}
|\nabla_{\mathrm{H}} e^{h \Delta_{\mathrm{H}}} f|(g) \le K e^{h \Delta_{\mathrm{H}}} \left( |\nabla_{\mathrm{H}} f| \right)(g), \quad \forall h > 0,  g \in \mathrm{H}, %{\color{red}\in G}, \
f \in C_0^{\infty}(\mathrm{H}),
\end{align}
where $\nabla_{\mathrm{H}}$ and $e^{h \Delta_{\mathrm{H}}}$ denote the horizontal gradient and the heat semigroup associated with the canonical sub-Laplacian $\Delta_{\mathrm{H}}$ of the H-type group $\mathrm{H}$ respectively,
are obtained in \cite{HL10} and \cite{E10} independently.

Moreover, in \cite{TM17}, Bruno and Calzi give asymptotic estimates at infinity of all radial partial
derivatives of the heat kernel on H-type groups.

The paper is organized as follows.  In Section \ref{st2} we prove the validity of differentiation of a general asymptotic expansion with respect to a particular kind of asymptotic sequence. A couple of examples are given in Section \ref{st3}. After a review of H-type groups in Section \ref{st4}, we give a direct application of the method of differentiation to a special case of asymptotics for the heat kernel on H-type groups in Section \ref{st5}.
For the general case, we present the theorem and make some preparations in Section \ref{st6}. The proof of the theorem will be carried out in
Section \ref{st7}.  In Section \ref{st8}, we give the proof of our main lemma, which completes the proof in Section \ref{st7}.

\section{Main results}\label{st2}
\setcounter{equation}{0}
Without loss of generality, we can assume $\Delta = \{z: |z| > R > 0, |\arg(z)|< \theta_0 < \pi\}$. We begin with the definition of
a class of analytic functions on $\Delta$.
\begin{definition}
For given $\theta_1 \in ( 0, \theta_0)$, let $\Delta^{\prime} =  \Delta  \cap \{z: |\arg(z)| \le \theta_1 \}$. A good test function on $\Delta^{\prime}$ is an analytic function $g$ on $\Delta$ which satisfies the following condition: there exist two constants $C > 0$ and
$R_1 > R$ such that $\forall z \in
\Delta^{\prime} \cap \{z : |z| \ge R_1\}$, there exists a constant $R(z) \in \left(0 , \frac{1}{2}|z|\right]$ satisfying
$C_z = \{\xi: |\xi - z| = R(z) \} \subset \Delta$ and
\begin{align}\label{E}
\int_{C_z}\frac{|g(\xi)|}{|\xi - z|^2}|d\xi| \le C |g^{\prime}(z)|.
\end{align}
The class of good test functions on $\Delta^{\prime}$ is denoted by $GTF(\theta_0,\theta_1,R)$.
\end{definition}

For the property of good test functions on $\Delta^{\prime}$, we introduce the following lemma.

\begin{lemma}\label{l1}
If $g \in GTF(\theta_0,\theta_1,R)$ and $f$ is an analytic function on $\Delta$ satisfying
\begin{align*}
f(z) \sim g(z),  \qquad \mbox{$z \to \infty$ in $\Delta$,}
\end{align*}
that is,
\begin{align*}
\lim_{z \to \infty, z \in \Delta} \frac{f(z)}{g(z)} = 1,
\end{align*}
then we have
\begin{align*}
f^{\prime}(z) \sim g^{\prime}(z),  \qquad \mbox{$z \to \infty$ in $\Delta^{\prime}$.}
\end{align*}
\end{lemma}

\noindent \textbf{Proof.} For any given $\epsilon > 0$, there is a constant $R_0 > R$ such that
\begin{align*}
|f(z) - g(z)| < \epsilon |g(z)|, \qquad \forall  z \in \Delta,  |z| \ge R_0.
\end{align*}
From our assumption, there are constants $C>0$ and $R_1 > R$ such that
$\forall z \in \Delta^{\prime} \cap \{z : |z| \ge R_1\}$, there is a constant $R(z) \in \left(0 , \frac{1}{2}|z|\right]$ satisfying
$C_z = \{\xi: |\xi - z| = R(z) \} \subset \Delta$ and \eqref{E}.
Let $R^{\prime} = \max(2R_0, R_1)$, and then for any $z \in \Delta^{\prime} \cap \{z: |z| \ge R^{\prime}\}$, we have $C_z \subset \Delta \cap
\{z: |z| \ge R_0\}$ by the choice of $R^{\prime}$ and $R(z)$. By Cauchy's integral formula, for any $z \in \Delta^{\prime} \cap \{z: |z| \ge R^{\prime}\}$,
\begin{align*}
|f^{\prime}(z) - g^{\prime}(z)|& =\left| \frac{1}{2\pi \imath} \int_{C_z} \frac{f(\xi) - g(\xi)}{(\xi - z)^2} d\xi \right| \\
&\le \frac{1}{2\pi}\int_{C_z} \frac{|f(\xi) - g(\xi)|}{|\xi - z|^2} |d\xi| \\
&\le \frac{\epsilon}{2\pi} \int_{C_z} \frac{|g(\xi)|}{|\xi - z|^2} |d\xi| \le \frac{C\epsilon}{2\pi}|g^{\prime}(z)|,
\end{align*}
which finishes our proof.
~ \hspace*{20pt} ~ \hfill $\Box$

In particular, we obtain the following corollary, which is crucial in the proof our main theorem.

\begin{corollary}\label{c1}
If $g \in GTF(\theta_0,\theta_1,R)$ and $f$ is an analytic function on $\Delta$ satisfying
\begin{align*}
f(z) = o(g(z)),  \qquad \mbox{$z \to \infty$ in $\Delta$,}
\end{align*}
then we have
\begin{align*}
f^{\prime}(z) = o(g^{\prime}(z)),  \qquad \mbox{$z \to \infty$ in $\Delta^{\prime}$.}
\end{align*}
\end{corollary}
\noindent \textbf{Proof.}
Notice that our assumption and result are equivalent to
\begin{align*}
f(z) + g(z) \sim g(z),  \qquad \mbox{$z \to \infty$ in $\Delta$,}
\end{align*}
and
\begin{align*}
f^{\prime}(z) + g^{\prime}(z) \sim g^{\prime}(z),  \qquad \mbox{$z \to \infty$ in $\Delta^{\prime}$,}
\end{align*}
respectively. Then the corollary follows from Lemma \ref{l1}.
~ \hspace*{20pt} ~ \hfill $\Box$

\vskip4pt

Analogous to the proof of Lemma \ref{l1}, we can also obtain a counterpart of Corollary \ref{c1}.

\begin{corollary}\label{c2}
If $g \in GTF(\theta_0,\theta_1,R)$ and $f$ is an analytic function on $\Delta$ satisfying
\begin{align*}
f(z) = O(g(z)),  \qquad \mbox{$z \to \infty$ in $\Delta$,}
\end{align*}
then we have
\begin{align*}
f^{\prime}(z) = O(g^{\prime}(z)),  \qquad \mbox{$z \to \infty$ in $\Delta^{\prime}$.}
\end{align*}
\end{corollary}

Before we state our main theorem, we introduce the definition of good asymptotic sequences.

\begin{definition}
For given $\theta_1 \in ( 0, \theta_0)$ and $\Delta^{\prime} =  \Delta  \cap \{z: |\arg(z)| \le \theta_1 \}$, a good asymptotic sequence
on $\Delta^{\prime}$ is an asymptotic sequence $\{g_n\}$ as $z \to \infty$ in $\Delta$ (that is, a sequence of functions such that for any $n \ge 0$ we have $g_{n+1} = o(g_n)$ as $z \to \infty$ in $\Delta$) such that for any $n$, $g_n \in GTF(\theta_0,\theta_1,R)$. In other words, each $g_n$ is
a good test function on  $\Delta^{\prime}$. The class of good asymptotic sequences on $\Delta^{\prime}$ is denoted by $GAS(\theta_0,\theta_1,R)$.
\end{definition}

The reason why we call an asymptotic sequence as $z \to \infty$ in $\Delta$ a good asymptotic sequence on $\Delta^{\prime}$ can be answered by the following theorem.

\begin{theorem}\label{t1}
If $\{g_n\} \in GAS(\theta_0,\theta_1,R)$, then $\{g_n^{\prime}\}$ is an asymptotic sequence
as $z \to \infty$ in $\Delta^{\prime}$.
Furthermore, if $f$ and ${f_n}$ are analytic in $\Delta$ such that
\begin{align*}
f(z) \sim \sum_{n=0}^{\infty} f_n(z), \qquad \{g_n\},\qquad \mbox{$z \to \infty$ in $\Delta$,}
\end{align*}
then
\begin{align*}
f^{\prime}(z) \sim \sum_{n=0}^{\infty} f_n^{\prime}(z), \qquad \{g_n^{\prime}\},\qquad \mbox{$z \to \infty$ in $\Delta^{\prime}$.}
\end{align*}
\end{theorem}

\noindent \textbf{Proof.}
This theorem follows directly from Corollary \ref{c1}. For example, for any $n$ we have
$g_{n+1}(z) = o(g_n(z))$ as $z \to \infty$ in $\Delta$. Then from Corollary \ref{c1} we have
$g_{n+1}^{\prime}(z) = o(g_n^{\prime}(z))$ as $z \to \infty$ in $\Delta^{\prime}$, which proves the first part of the theorem.
The second part of the theorem follows similarly.
~ \hspace*{20pt} ~ \hfill $\Box$

\section{Examples}\label{st3}
\setcounter{equation}{0}
So far we have worked in the abstract frame and we give some examples in this section. We first fix some notation.
\begin{definition}
Suppose $f$ and $g$ are two real-valued functions. The notation $f \lesssim g$ means that there exists a constant $C > 0$ such that
$f \leq C g$. Furthermore, we use the notation $f \approx g$ if $f \lesssim g$ and $g \lesssim f$.
\end{definition}

We give a sufficient condition which ensures a analytic function on $\Delta$ to be a good test function on $\Delta^{\prime}$.

\begin{lemma}\label{l2}
For given $\theta_1 \in ( 0, \theta_0)$ and $\Delta^{\prime} =  \Delta  \cap \{z: |\arg(z)| \le \theta_1 \}$, if an analytic function $g$ on $\Delta$ satisfies the following condition: there exist two constants $C > 0$ and $R_1 > R$ such that $\forall z \in
\Delta^{\prime} \cap \{z : |z| \ge R_1\}$, there exists a constant $R(z) \in \left(0 , \frac{1}{2}|z|\right]$ satisfying
$C_z = \{\xi: |\xi - z| = R(z) \} \subset \Delta$ and
\begin{align}
\frac{|g(\xi)|}{R(z)} \le \frac{C}{2\pi} |g^{\prime}(z)|, \qquad \forall \xi \in C_z,
\end{align}
then $g \in GTF(\theta_0,\theta_1,R)$.
Conversely, if $g \in GTF(\theta_0,\theta_1,R)$, there exist two constants $C > 0$ and $R_1 > R$ such that $\forall z \in
\Delta^{\prime} \cap \{z : |z| \ge R_1\}$ we have
\begin{align}
|g(z)| \le \frac{C}{4\pi} |zg^{\prime}(z)|.
\end{align}
\end{lemma}

\begin{example}\label{e1}
For any $0 <\theta_1 <\theta_0<\pi$, $R > 0$, $\alpha \in \C \backslash\{0\}$ and $\beta \in \C$ fixed, we have $z^{\alpha}\log^{\beta}{z} \in GTF(\theta_0,\theta_1,R)$.
\end{example}
\noindent \textbf{Proof.}
Set $R_1 = \max(3R,e^{10})$ and $R(z) = \frac{1}{2}\sin(\theta_0 - \theta_1)|z|$ for $z \in \Delta^{\prime} \cap \{z: |z| \ge R_1\}$.
Then for any $z \in \Delta^{\prime} \cap \{z: |z| \ge R_1\}$ and $\xi \in C_z$ we have $C_z \subset \Delta$ and $R(z) \approx |z| \approx |\xi|$.
As a result,  for any $z \in \Delta^{\prime} \cap \{z: |z| \ge R_1\}$ and $\xi \in C_z$, we have
\begin{align*}
|g(\xi)| \approx |\xi|^{\Re\alpha} |\log{\xi}|^{\Re \beta} \approx |z|^{\Re\alpha} |\log{z}|^{\Re \beta} \approx R(z)
\left|(z^{\alpha}\log^{\beta}{z})^{\prime}\right|,
\end{align*}
where $\Re w$ denotes the real part of a complex number $w$.
Then the assertion follows from the first part of Lemma \ref{l2}.
~ \hspace*{20pt} ~ \hfill $\Box$

\begin{remark}
From this example we deduce that for all $\theta_0$, $\theta_1$ and $R$ we have $\{z^{-n}\} \in GAS(\theta_0,\theta_1,R)$. Then Theorem \ref{t1}
is a generalization of the classical result.
\end{remark}

\begin{example}\label{e2}
For any $0 < \theta_1 < \theta_0 < \pi$, $R > 0$, $\alpha \in \C$, $\beta \in \C \backslash \{0\}$ and $\Re \gamma > 0$ fixed, we have $z^{\alpha}e^{\beta z^{\gamma}} \in GTF(\theta_0,\theta_1,R)$.
\end{example}
\noindent \textbf{Proof.}
We first pick $R_1 \ge \max(3R,e^{10}) >0$ large enough such that $\frac{1}{2}\sin(\theta_1 - \theta_0) R_1^{\Re \gamma} \ge 1$.
For $z \in \Delta^{\prime} \cap \{z: |z| \ge R_1\}$, we set $R(z) = |z|^{1 - \Re \gamma } \le \frac{1}{2}|z|\sin(\theta_0 - \theta_1)$.
Then for any $z \in \Delta^{\prime} \cap \{z: |z| \ge R_1\}$ and $\xi \in C_z$, we have $C_z \subset \Delta$ and $|\xi| \approx |z|$.
Observe that for any $z \in \Delta^{\prime} \cap \{z: |z| \ge R_1\}$ and $\xi \in C_z$,
\begin{align*}
\xi^{\gamma} = (z + (\xi - z))^{\gamma} = z^{\gamma} \left( 1 + \frac{\xi - z}{z}\right)^{\gamma}
= z^{\gamma}\left(1 + O\left(\frac{\xi - z}{z}\right)\right) = z^{\gamma} + O(1).
\end{align*}
As a result,
\begin{align*}
|g(\xi)| \approx |\xi|^{\Re \alpha} |e^{\beta z^{\gamma}}| |e^{O(1)}| \approx |z|^{\Re \alpha} |e^{\beta z^{\gamma}}|
\approx R(z) \left|(z^{\alpha}e^{\beta z^{\gamma}})^{\prime}\right|.
\end{align*}
Then the assertion follows from the first part of Lemma \ref{l2}.
~ \hspace*{20pt} ~ \hfill $\Box$

\begin{example}\label{e3}
For any $\theta_0$, $\theta_1$ and $R$, $\log{z}  \notin GTF(\theta_0,\theta_1,R)$.
\end{example}
\noindent \textbf{Proof.}
If for some $\theta_0$, $\theta_1$ and $R$, $\log{z}  \in GTF(\theta_0,\theta_1,R)$, then from the second part of Lemma \ref{l2} we have
$|\log{z}| \lesssim 1$ for $z$ real and large enough, which is a contradiction.
~ \hspace*{20pt} ~ \hfill $\Box$

\begin{remark}
Furthermore, Lemma \ref{l1} is not valid for $\log{z}$, that is,
if $f(z) \sim \log{z}$ as $z \to \infty$ in $\Delta$ for some analytic function $f$, then $f^{\prime}(z) \nsim \frac{1}{z}$ as $z \to \infty$ in $\Delta^{\prime}$.
For example, $\log{z} + \sin(\log{z}) \sim \log{z}$
as $z \to \infty$ in $\Delta$ but $\frac{1}{z} + \frac{\cos(\log{z})}{z} \sim \frac{1}{z}$ fails even on the real line.
However, if $f(z) \sim \log{z}$ as $z \to \infty$ in $\Delta$, then $zf(z) \sim z\log{z}$ as $z \to \infty$ in $\Delta$.
Since $z\log{z} \in GTF(\theta_0,\theta_1,R)$ by Example \ref{e1}, we have $f(z) + zf^{\prime}(z) \sim \log{z}$ as $z \to \infty$ in $\Delta^{\prime}$,
which yields $f^{\prime}(z) = O\left(\frac{\log{z}}{z}\right)$ as $z \to \infty$ in $\Delta^{\prime}$.
\end{remark}

To apply Lemma \ref{l2} we can use the following lemma.

\begin{lemma}\label{l3}
Assume $\Sigma \subset \C$ is a nonempty convex open set and $g$ is an analytic function on $\Sigma$ satisfying the following condition:
For some constant $0 < C_1 < C_2$, we have
\begin{align*}
C_1 \le \left|\frac{g^{\prime}(z)}{g(z)} \right| \le C_2.
\end{align*}
Then we have
\begin{align*}
|g(\xi)| \le C_1^{-1} e^{C_2|\xi - z|} |g^{\prime}(z)|, \qquad \forall \xi,z \in \Sigma.
\end{align*}
\end{lemma}

\noindent \textbf{Proof.} Pick a $z_0 \in \Sigma$, we define
\begin{align*}
\phi(z) = \int_{z_0}^z \frac{g^{\prime}(s)}{g(s)} ds, \qquad \forall z \in \Sigma,
\end{align*}
where the integral is along any curve linking $z_0$ and $z$. It is well-defined and analytic since $\frac{g^{\prime}(z)}{g(z)}$ is analytic
on $\Sigma$ from our assumption.

As a result, we have $g^{\prime}(z) = \phi^{\prime}(z) g(z)$, which is equivalent to
\begin{align*}
\frac{d}{dz}(e^{-\phi(z)}g(z)) = 0
\end{align*}
and finally $g(z) = Ke^{\phi(z)}$ for some constant $K$.

Then we have
\begin{align*}
|e^{\phi(\xi)}| &=|e^{\phi(z)}||e^{\phi(\xi) - \phi(z)}| \le
|e^{\phi(z)}|e^{|\phi(\xi) - \phi(z)|}\\
 &\le |e^{\phi(z)}|e^{C_2|\xi - z|} \le e^{C_2|\xi - z|} C_1^{-1} |\phi^{\prime}(z)e^{\phi(z)}|,
\end{align*}
which proves the lemma.
~ \hspace*{20pt} ~ \hfill $\Box$

\begin{remark}
From Lemma \ref{l3}, for any $0 < \theta_1 < \theta_0 < \pi$, $R > 0$, $\alpha \in \C$ and $\beta \in \C \backslash \{0\}$ fixed we have $z^{\alpha}e^{\beta z} \in GTF(\theta_0,\theta_1,R)$, which is a special case of Example \ref{e2}.
\end{remark}

\section{Preliminaries for H-type groups}\label{st4}
\setcounter{equation}{0}

\subsection{H-type groups}
Recall that an H-type group can be considered as $\H(2n,m) = \R^{2n} \times \R^{m} (m,n \in \N)$ with the group law (see Theorem A.2, p. 199 of \cite{BU04};
see also \cite{K80} for an original definition)
\begin{align*}
(z,t)\cdot(z^{\prime},t^{\prime}) = (z + z^{\prime}, t+t^{\prime} + 2^{-1}\langle z,Uz^{\prime} \rangle)
\end{align*}
with $z = (z_1, \ldots, z_{2n}) \in \R^{2n}$, $t = (t_1, \ldots, t_m) \in \R^m$ and
\begin{align*}
\langle z,Uz^{\prime} \rangle = (\langle z,U^{(1)}z^{\prime} \rangle, \ldots, \langle z,U^{(m)}z^{\prime} \rangle) \in \R^m,
\end{align*}
where the matrices $U^{(1)}, \ldots, U^{(m)}$ have the following properties:

1. $U^{(j)} (1 \le j \le m)$ is a $(2n)\times(2n)$ skew-symmetric and orthogonal matrix.

2. $U^{(i)}U^{(j)} + U^{(j)}U^{(i)} = 0$ for all $1 \le i \ne j \le m$.

Let $U^{(j)} = (U^{(j)}_{k,l})_{1 \le k,l \le 2n}(1 \le j \le m)$. When $m = 1$, $\H(2n,1)$ is often called the (isotropic) Heisenberg group of real dimension $2n + 1$. The canonical sub-Laplacian on $\H(2n,m)$ is given by
sum of squares: $\Delta_{\H(2n,m)} = \sum\limits_{l=1}^{2n}\mathrm{X}_l^2$, where $\mathrm{X}_l(1 \le l \le 2n)$ are the left-invariant vector fields on $\H(2n,m)$,
defined by
\begin{align*}
\mathrm{X}_l = \frac{\partial}{\partial z_l} + \frac{1}{2} \sum_{j=1}^m \left(\sum_{k=1}^{2n} z_k U^{(j)}_{k,l}\right)\frac{\partial}{\partial t_j}.
\end{align*}
Let $\nabla_{\H(2n,m)} = (\mathrm{X}_1, \ldots, \mathrm{X}_{2n})$ be the horizontal gradient, $d$ the Carnot-Carath\'{e}odory distance, and $p_h(h > 0)$ the
heat kernel on $\H(2n,m)$. We let $o = (0,0)$ be the identity element of $\H(2n,m)$; $g = (z,t) \in \R^{2n} \times \R^m$ denotes a point in $\H(2n,m)$.

For simplicity we also denote
\begin{align*}
d(g) = d(g,o) \quad \mbox{and} \quad p_h(g) = p_h(g,o).
\end{align*}

\subsection{The Carnot-Carath\'{e}odory distance}
Let $\mu: (-\pi, \pi) \to (-\infty, \infty)$ defined by
\begin{align*}
\mu(w) = \frac{2w - \sin(2w)}{2\sin^2{w}},
\end{align*}
which is a strictly increasing diffeomorphism.

It is well-known that(see \cite{Li10} or \cite{HL10} for example)
\begin{align*}
d^2(z,t)=\left\{\begin{array}{ll}
\left(\frac{\theta}{\sin{\theta}}\right)^2|z|^2 & \textrm{if $z \ne 0$, where $\theta=\mu^{-1}\left(\frac{4|t|}{|z|^2}\right)$}\\
4\pi|t| & \textrm{if $z=0$}
\end{array} .\right.
\end{align*}

From the expression above we have
\begin{eqnarray*}
d\left(\frac{z}{\sqrt{h}},\frac{t}{h}\right)=\frac{1}{\sqrt{h}}d(z,t), \quad \forall h > 0.
\end{eqnarray*}

Moreover, using the equivalence between the Carnot-Carath\'eodory distance and a homogeneous norm on stratified groups (see for example \cite{VSC92}), or by a
direct calculation, we have
\begin{align*}
d^2(z,t) \approx |z|^2 + |t|.
\end{align*}

\subsection{Heat kernel on H-type groups}
For $u,v \ge 0$, we define
\begin{align}\label{defp}
p(n,m;u,v) = \frac{2}{(4\pi)^{n+\frac{m}{2}}} \int_0^{\infty} s^{m-1}
e^{-\frac{u}{4}s\coth{s}} \left(\frac{s}{\sinh{s}}\right)^n \left(\frac{sv}{2}\right)^{-\frac{m-2}{2}} J_{\frac{m-2}{2}}(sv) ds.
\end{align}

It is well known that (see \cite{Li10} for example)
\begin{align*}
p_h(z,t)
= \frac{1}{h^{n+m}} p\left(n,m;\frac{|z|^2}{h},\frac{|t|}{h}\right).
\end{align*}

As a result, to obtain asymptotic estimates at infinity of the heat kernel on H-type groups, it is equivalent to consider asymptotic estimates at infinity of the function
$p(n,m;\cdot,\cdot)$.

Here are some useful facts about the function $p(n,m;\cdot,\cdot)$, which can be found in \cite[(1.12),(1,13)]{Li10}:
\begin{align}\label{key1}
\frac{\partial}{\partial v}p(n,m;u,v) &= -2\pi v p(n,m+2;u,v),\\
\label{key2}
p(n,m;u,v) &= 2 \int_v^{\infty} \frac{h}{\sqrt{h^2 - v^2}}p(n,m+1;u,h) dh.
\end{align}

\section{A direct application to the case $u = 0$} \label{st5}
\setcounter{equation}{0}
We can extend our result to other regions in the complex plane $\C$. For example, if $\Omega = \{z: \Re z > \log{R}, |\Im z| < \theta_0 < \pi \}$, where $\Im w$ denotes the
imaginary part of a complex number $w$. Then for analytic functions $f$ and $g$ on $\Omega$ we have $f(z) \sim g(z)$ as $z \to \infty$ in $\Omega$ if and only if $f(\log{z}) \sim g(\log{z})$
as $z \to \infty$ in $\Delta$. Moreover, if $g(\log{z}) \in GTF(\theta_0,\theta_1,R)$, then we have
$\frac{f^{\prime}(\log{z})}{z} \sim \frac{g^{\prime}(\log{z})}{z}$ as $z \to \infty$ in $\Delta^{\prime}$, which yields
$f^{\prime}(z) \sim g^{\prime}(z)$ as $z \to \infty$ in $\Omega^{\prime} = \{z: \Re z > \log{R}, |\Im z| \le \theta_1 < \theta_0 \}$.
In particular, if $\alpha \ne 0$, $g(z) = z^{\beta}e^{\alpha z}$ has the property that for any analytic function $f$ on $\Omega$ such that
$f \sim g$ as $z \to \infty$ in $\Omega$, then $f^{\prime} \sim g^{\prime}$ as $z \to \infty$ in $\Omega^{\prime}$ (see Example \ref{e1}). Other regions like
$\Theta = \{z: |\arg{z}| < \theta_0 < \pi, 0<|z|<L\}$ can be handled similarly (with $z \to \infty$ replaced by $z \to 0$).

Another approach is to adapt the definition of good test functions on $\Delta^{\prime}$ to other regions.
For example, we can define a good test function on $\Omega^{\prime}$ by an analytic function $f$ on $\Omega$ with the following property:
there exist constants $C > 0$ and $R_1 > R$ such that $\forall z \in \Omega^{\prime} \cap \{z: \Re z > \log{R_1} \}$, there exists a constant
$R(z) \in (0, \frac{1}{2}(\theta_0 - \theta_1)]$ satisfying $C_z = \{\xi: |\xi - z| = R(z) \} \subset \Omega$ and \eqref{E}.
Then we can deduce counterparts of Lemma \ref{l1} and Theorem \ref{t1} without difficulties. \\

Let $b(n,m;v) = p(n,m;0,v)$ where $p$ is defined in \eqref{defp}. From the well-known result $J_{-\frac{1}{2}}(z)=
\sqrt{\frac{2}{\pi z}}\cos(z)$ (cf. \cite[8.464.2, p. 924]{GR07}) we have
\begin{align*}
b(n,1;v)=\frac{1}{(4\pi)^{n+\frac{1}{2}}\sqrt{\pi}}\int_{\R} \left(\frac{s}{\sinh{s}}\right)^n e^{\imath v s } ds,
\end{align*}
which is analytic on the strip $\{v: |\Im v|< \delta \}$ for sufficient small $\delta > 0$ by the property of Laplace's transform
(cf. \cite[Theorem 5a, p. 57]{W41}).
By Cauchy's theorem,
\begin{align*}
b(n,1;v) = c_n \int_{\R + \imath \frac{3\pi}{2}} \left(\frac{s}{\sinh{s}}\right)^n e^{\imath v s } ds + c_n 2\pi \imath \mathrm{Res}
\left[ \left(\frac{s}{\sinh{s}}\right)^n e^{\imath v s }, \imath \pi\right] = M_1 + M_2,
\end{align*}
where $c_n = \frac{1}{(4\pi)^{n+\frac{1}{2}}\sqrt{\pi}}$.
For the integral, notice that $\sinh(\zeta + \imath \eta) = \sinh{\zeta}\cos{\eta} + \imath \cosh{\zeta}\sin{\eta}$, which yields
\begin{align*}
|M_1| \lesssim e^{-\Re{v} \frac{3\pi}{2}} \int_{\R} \frac{(|s|+1)^n}{(\cosh{s})^n} e^{\delta |s|} ds \lesssim e^{-\Re{v} \frac{3\pi}{2}} = o(v^{n-1}e^{-\pi v}),
\qquad v \to +\infty \mbox{ in } \{v: |\Im v|< \delta \}.
\end{align*}
For the residue, a direct computation gives
\begin{align*}
2\pi \imath \mathrm{Res}\left[ \left(\frac{s}{\sinh{s}}\right)^n e^{\imath v s }, \imath \pi\right]
= \frac{2\pi^{n+1}}{(n-1)!}v^{n-1}e^{- \pi v}(1  +  o(1)), \qquad v \to +\infty \mbox{ in } \{v: |\Im v|< \delta \}.
\end{align*}
After all, $b(n,1;v) \sim \frac{1}{4^n (n - 1)!} v^{n-1}e^{- \pi v}$ as $v \to \infty$ in the strip $\{v: \Re v > 0, |\Im v|< \delta \}$.
From \eqref{key1}, we obtain
\begin{align*}
-2\pi v b(n,m+2;v) = \frac{d}{d v} b(n,m;v),
\end{align*}
which yields
\begin{align*}
b(n,2k+1;v) \sim \frac{1}{2^k 4^n (n - 1)!} v^{n-k-1} e^{- \pi v}, \qquad v \to +\infty \mbox{ in } \left\{v: |\Im v|< \delta_k = \frac{2^k \delta}{3^k}\right\}
\end{align*}
from the property of $z^{\beta}e^{\alpha z}$ when $\alpha \ne 0$ discussed above.
Finally, from \eqref{key2}, we obtain
\begin{align*}
b(n,2k;v) &\sim \frac{2}{2^k 4^n (n - 1)!} \int_v^{\infty} \frac{s^{n-k}}{\sqrt{s^2 - v^2}} e^{- \pi s} ds  \\
&= \frac{2}{2^k 4^n (n - 1)!} v^{n - k}e^{- \pi v} \int_0^{\infty} \frac{(r + 1)^{n - k}}{\sqrt{r} \sqrt{r+2}} e^{- \pi v r} dr \\
&=  \frac{\sqrt{2}}{2^k 4^n (n - 1)!} v^{n - k}e^{- \pi v} \left( \int_0^1 \frac{1 + O(r)}{\sqrt{r}} e^{- \pi v r} dr
+ o\left(v^{-\frac{1}{2}}\right) \right) \\
&= \frac{\sqrt{2}}{2^k 4^n (n - 1)!} v^{n - k - \frac{1}{2}}e^{- \pi v}(1 + o(1))
\end{align*}
as $v \to +\infty$ on the real line.

\section{Complexification and asymptotics in the general case $u > 0$} \label{st6}
\setcounter{equation}{0}
Before we give the asymptotic estimate at infinity of the heat kernel on H-type groups when $u > 0$, we need some
discussion about the analyticity.
Recall that $\mu: (-\pi, \pi) \to (-\infty, \infty)$ is defined by
\begin{align*}
\mu(w) = \frac{2w - \sin(2w)}{2\sin^2{w}}.
\end{align*}
However, we can regard $\mu$ as a meromorphic function on $\C$.
Notice that when $w$ is near $\pi$, we have
\begin{align}\label{mu}
\mu(w) = \frac{\pi}{(\pi - w)^2 } (1 + o(1)).
\end{align}

As a consequence, picking a suitable branch, $\mu^{-1}$ is defined and is an analytic function from
$\Omega_1 = \{z| |z|>R_0>100, |\arg{z}|< \eta_0< \frac{\pi}{4}\}$ to $\Omega_2 = \{w| |\pi - w|<R_0^{\prime}<1, |\arg(\pi - w)|< \eta_0^{\prime}
<\frac{\pi}{8}\}$.

Now we let $v$ be a complex number. Then when $u>0$, $\Re v \gg 1$, $\Re v \gg u$ and $|\Im v| \ll 1$, we have
$\frac{4v}{u} \in \Omega_1$ and we define $\theta = \mu^{-1}\left(\frac{4v}{u}\right) \in \Omega_2$.
We remark that $\theta$ is analytic in $v$, $\frac{d\theta}{dv} = \frac{4}{u \mu^{\prime}(\theta)}$ and $\theta \to \pi$ under our assumption.
Furthermore, we define $d^2(u,v) = \left(\frac{\theta}{\sin\theta}\right)^2 u $, which is again an analytic function in $v$ and
$\frac{d}{dv}d^2(u,v) = \theta \mu^{\prime}(\theta) u \frac{d\theta}{dv} = 4 \theta$.
Let $\epsilon = \pi - \theta$. Then we have $\epsilon \to 0$, $|\arg{\epsilon}| < \eta_0^{\prime}$ and
\begin{align*}
\epsilon = \sqrt{\frac{\pi u }{4 v}} (1 + o(1))
\end{align*}
from (\ref{mu}).

\vskip4pt

The following theorem is our main result of the asymptotic behaviour at infinity of the heat kernel of H-type groups in this paper.

\begin{theorem}\label{t2}
For $u > 0$, $v \gg 1$ and $v \gg u$, we have
\begin{align}\label{asyp}
p(n,m;u,v) = 2^{-2n - \frac{m-1}{2}} v^{-\frac{m-1}{2}} e^{-\frac{d^2(u,v)}{4}} \epsilon^{1-n} e^{-\frac{\pi}{2\epsilon}u} I_{n-1}\left(\frac{\pi}{2\epsilon}u\right)(1 + o(1)),
\end{align}
where the modified Bessel function $I_{\nu}$ ($\nu > -\frac{1}{2}$) is defined as
\begin{eqnarray*}
I_{\nu}(u)  =  \frac{1}{\sqrt{\pi} \Gamma(\nu + \frac{1}{2})} \Big(
\frac{u}{2} \Big)^{\nu} \int_{-1}^1 (1 - h^2)^{\nu - \frac{1}{2}}
e^{-u h} \, dh, \quad u \in \C.
\end{eqnarray*}
\end{theorem}

\begin{remark}
The $o(1)$ in \eqref{asyp} only depends on $m$, $n$, $v$ and $\frac{v}{u}$, not directly on $u$.
\end{remark}

\begin{remark}
The remaining case (the case when $u > 0$, $u \gg 1$ and $\frac{v}{u} \lesssim 1$) can be solved by the method of
stationary phase. For more details, one can refer to \cite{Li07} or \cite{E09}.
\end{remark}

Although this theorem only considers the behaviour of $v$ when $v$ is real, to apply our method of differentiation we need $v$ to be complex.
However, when $m=1$, that is, in the isotropic Heisenberg group case, a direct modification of \cite{Li07} gives the next lemma (which will be
proved in the Appendix for the sake of completeness).
\begin{lemma}\label{l4}
There exists a constant $\delta \in (0,1)$ such that for $u > 0$, $\Re v \gg 1$, $\Re v \gg u$ and $|\Im v| < \delta$ we have
\begin{align}\label{asyp1}
p(n,1;u,v) = 2^{-2n}  e^{-\frac{d^2(u,v)}{4}} \epsilon^{1-n} e^{-\frac{\pi}{2\epsilon}u} I_{n-1}\left(\frac{\pi}{2\epsilon}u\right)(1 + o(1)),
\end{align}
where $o(1)$ in \eqref{asyp1} only depends on $n$, $\delta$, $\Re v$ and $\Re \frac{v}{u}$.
\end{lemma}

\section{Proof of Theorem \ref{t2}} \label{st7}
\setcounter{equation}{0}
We first prove the theorem when $m$ is odd by induction. When $m = 1$, it is the result of Lemma \ref{l4}.
We assume that for some odd $m = 2k + 1$, $u > 0$, $\Re v \gg 1$, $\Re v \gg u$ and $|\Im v| < \delta_k < 1$ we have
\begin{align*}
p(n,2k+1;u,v) = 2^{-2n - k} v^{-k} e^{-\frac{d^2(u,v)}{4}} \epsilon^{1-n} e^{-\frac{\pi}{2\epsilon}u} I_{n-1}\left(\frac{\pi}{2\epsilon}u\right)(1 + o(1))
\end{align*}
and we need to prove it is also true for $p(n,2k+3;u,v)$ (with $\delta_k$ replaced by $\delta_{k+1} = \frac{2}{3}\delta_{k}$).
We denote the function on the right hand side by $q(n,2k+1;u,v)$. Then from our discussion above, for fixed $u > 0$, $q$ is an analytic function in $v$. Since $p(n,1;u,v)$ is analytic in $v$ on $\{v: |\Im v| < \delta\}$ (see Appendix), from \eqref{key1}, we have
$p(n, 2k+1;u,v)$ is analytic in $v$ on  $\{v: \Re v > 0, |\Im v| < \delta\}$.
%However, by the property of the Laplace transform, for fixed $u > 0$, $p(n,m;u,v)$ is an analytic function in $v$ when $|\Im v| < \delta_0$.
Recall that (cf. \cite[\S 8.445, p. 919]{GR07} and \cite[\S8.451.5, p. 920]{GR07}) for fixed $\nu > 0$,
\begin{align}\label{I1}
I_{\nu -1}(z) &= \frac{1}{\Gamma(\nu)} \left(\frac{z}{2}\right)^{\nu - 1} (1 + o(1)), \qquad z \to 0, \\
\label{I2}
I_{\nu -1}(z) &= e^{z} \frac{1}{\sqrt{2 \pi z}}(1 + o(1)),\qquad z \to \infty,
\end{align}
and
\begin{align}\label{key3}
d^2(u,v) = 4 \pi v (1 + o(1)) = \frac{\pi^2}{\epsilon^2} u (1 + o(1)), \qquad \Re v \to +\infty, \Re\frac{v}{u} \to +\infty.
\end{align}

Considering separately the case $|\frac{u}{\epsilon}| < 1$ and $|\frac{u}{\epsilon}| \ge 1$ and using \eqref{I1}--\eqref{key3}
we have
\begin{align}\label{propq}
\frac{\partial}{\partial v} q(n,2k+1;u,v) &= \left( -\theta
-  \frac{k}{v} - \frac{2\pi}{\epsilon^2 \mu^{\prime}(\theta)} +  \frac{2\pi}{\epsilon^2 \mu^{\prime}(\theta)}
\frac{I_{n-2}\left(\frac{\pi}{2\epsilon}u\right)}{I_{n-1}\left(\frac{\pi}{2\epsilon}u\right)} \right)q(n,2k+1;u,v) \\
&= \left( -\theta + O(|t|^{-1}) + O(|\epsilon|)\right) q(n,2k+1;u,v).
\end{align}
Then for any $\rho > 0$, from assumption we get $R_1 \ge 1$ and $R_2 \ge 1$ such that
\begin{align*}
\left|\frac{p(n,2k+1;u,v)- q(n,2k+1;u,v)}{q(n,2k+1;u,v)}\right| \le \rho
\end{align*}
and
\begin{align*}
\frac{\pi}{2} \le \left| \frac{\partial_v q(n,2k+1;u,v)}{q(n,2k+1;u,v)} \right| \le 2 \pi
\end{align*}
when $\Re v > R_1$, $\Re v > R_2 u$ and $|\Im v|< \delta_k <1 $.
Then when $\Re v > 2R_1$, $\Re v > 2R_2 u$, $|\Im v|< \frac{2\delta_k}{3} = \delta_{k+1}  $ and $|\xi - v| = \frac{\delta_k}{3}$ we have
$\Re \xi > R_1$, $\Re \xi > R_2 u$, $|\Im \xi|< \delta_k <1 $ and
\begin{align*}
|q(n,2k+1;u,\xi)| \le \frac{2}{\pi} e^{\frac{2\pi\delta_k}{3}}|\partial_v q(n,2k+1;u,v)|
\end{align*}
from Lemma \ref{l3}.

As a result, we have
\begin{align*}
|\partial_v p(n,2k+1;u,v) - \partial_v q(n,2k+1;u,v)| &\le \frac{1}{2\pi} \int_{|\xi - v| = \frac{\delta_k}{3}} \frac{|p(n,2k+1;u,\xi) - q(n,2k+1;u,\xi)|}{|\xi - v|^2} |d\xi|\\
&\le \frac{\rho}{2\pi} \int_{|\xi - v| = \frac{\delta_k}{3}} \frac{|q(n,2k+1;u,\xi)|}{|\xi - v|^2} |d\xi| \\
&\le \frac{6\rho}{\pi \delta_k} e^{\frac{2\pi\delta_k}{3}}|\partial_v q(n,2k+1;u,v)|,
\end{align*}
which proves our result by (\ref{key1}).

For $m = 2k$ even, notice that from \eqref{propq} we have
\begin{align*}
\partial_v \log\left(q(n,2k+1;u,v)\right) = -(\pi + o(1)), \qquad v \to +\infty, \frac{v}{u} \to +\infty,
\end{align*}
which yields
\begin{align}\label{propq2}
\frac{q(n,2k+1;u,s)}{q(n,2k+1;u,v)} = e^{-(\pi + o(1))(s - v)}, \quad s \ge v,  v \to +\infty, \frac{v}{u} \to +\infty.
\end{align}

From \eqref{key2} and \eqref{propq2}, we obtain
\begin{align*}
p(n,2k;u,v) &\sim 2 \int_v^{\infty} \frac{s}{\sqrt{s^2 - v^2}} q(n,2k+1;u,s) ds  \\
&= 2q(n,2k+1;u,v) \int_v^{\infty} \frac{s}{\sqrt{s^2 - v^2}} e^{-(\pi + o(1))(s - v)} ds\\
&= 2vq(n,2k+1;u,v) \int_0^{\infty} \frac{(r + 1)}{\sqrt{r} \sqrt{r+2}} e^{- (\pi +o(1)) v r} dr \\
&= \sqrt{2} v q(n,2k+1;u,v) \left( \int_0^1 \frac{1 + O(r)}{\sqrt{r}} e^{- (\pi + o(1)) v r} dr
+ o\left(v^{-\frac{1}{2}}\right) \right) \\
&= \sqrt{2v}  q(n,2k+1;u,v)(1 + o(1)),
\end{align*}
which ends the proof of Theorem \ref{t2}.

\section{Appendix: proof of Lemma \ref{l4}} \label{st8}
As before, from the well-known result  $J_{-\frac{1}{2}}(z)=
\sqrt{\frac{2}{\pi z}}\cos(z)$  (cf. \cite[8.464.2, p. 924]{GR07}) we have
\begin{align}\label{p}
p(n,1;u,v)=\frac{1}{(4\pi)^{n+\frac{1}{2}}\sqrt{\pi}}\int_{\R} \left(\frac{s}{\sinh{s}}\right)^n e^{\imath v s -\frac{u}{4}s\coth{s}} ds.
\end{align}
Noticing that $\frac{u}{4}s\coth{s} \ge 0$, from the property of Laplace's transform (cf. \cite[Theorem 5a, p. 57]{W41}) there exists a constant $\delta$ (independent of $u$) such that
$p(n,1;u,v)$ is analytic on the strip $\{v: |\Im v|< \delta \}$. Moreover, we can choose $\delta$ small enough such that
when $u>0$, $\Re v \gg 1$ and $\Re v \gg u$, we have
$\frac{4v}{u} \in \Omega_1 = \{z| |z|>R_0>100, |\arg{z}|< \eta_0< \frac{\pi}{4}\}.$

We define
\begin{align*}
\phi(u,v;s) = \imath v s -\frac{u}{4}s\coth{s}.
\end{align*}
Then from direct computation we have
\begin{align*}
\phi(u,v;\imath \theta) = - \frac{d^2(u,v)}{4}, \qquad \frac{d}{ds}\phi(u,v;\imath \theta) = 0.
\end{align*}
It follows that
\begin{align}\label{phi}
\phi(u,v;\imath (\pi - \xi) ) = - (\pi - \xi) v + \frac{\pi u}{4 \xi} + G(u; \xi),
\end{align}
where
\begin{align*}
G(u; \xi) = \frac{u}{4} \left( \pi\left( \cot{\xi} - \frac{1}{\xi} \right)  - \xi \cot{\xi} \right)
\end{align*}
is an analytic function near the origin. We deduce from \eqref{phi} and $\varphi^{\prime}(u,v;\imath \theta) = 0$ that
\begin{align*}
v = \frac{\pi u}{4 \epsilon^2} - \frac{d}{d\xi} G(u;\epsilon).
\end{align*}
Combining this with \eqref{phi} we get
\begin{align*}
\phi(u,v;\imath (\pi - \xi)) - \phi(u,v;\imath (\pi - \epsilon)) =
\frac{\pi u}{4\epsilon} \left(\frac{\xi}{\epsilon} + \frac{\epsilon}{\xi} -2\right)
+ R(u,\epsilon;\xi),
\end{align*}
where
\begin{align}\label{R}
R(u,\epsilon;\xi) = G(u;\xi) - G(u; \epsilon) - \frac{d}{d\xi} G(u; \epsilon) (\xi - \epsilon).
\end{align}
From Cauchy's theorem we have
\begin{align*}
p(n,m;u,v)&=\frac{1}{(4\pi)^{n+\frac{1}{2}}\sqrt{\pi}}\left(\int_{\R + \imath \frac{3\pi}{2}} + \int_{|s - \imath \pi| = |\epsilon|} \left(\frac{s}{\sinh{s}}\right)^n e^{\imath v s -\frac{u}{4}s\coth{s}} ds \right)\\
&=\frac{1}{(4\pi)^{n+\frac{1}{2}}\sqrt{\pi}}(\mathrm{I} + \mathrm{II}) .
\end{align*}
We write $s = \zeta + \imath \eta$ with $\zeta, \eta \in \R$ and then
\begin{align*}
\Re(s \coth{s}) = \frac{\zeta \sinh{\zeta} \cosh{\zeta} + \eta \sin{\eta} \cos{\eta}}{ \sinh^2{\zeta} + \sin^2{\eta}}.
\end{align*}
It yields
\begin{align}\label{I}
|\mathrm{I}|  \lesssim e^{-\Re{v} \frac{3\pi}{2}} \int_{\R} \frac{(|s|+1)^n}{(\cosh{s})^n} e^{\delta |s|} ds \lesssim e^{-\Re{v} \frac{3\pi}{2}}.
\end{align}
Considering separately the case $|\frac{u}{\epsilon}| < 1$ and $|\frac{u}{\epsilon}| \ge 1$ and using \eqref{I1}--\eqref{key3}
we have
\begin{align}\label{R1}
\mathrm{I} = o\left(v^{-\frac{1}{2}}e^{-\frac{d^2(u,v)}{4}}\right) = o(q(n,1;u,v)).
\end{align}
Now we are in a position to estimate $\mathrm{II}$, and the proof is divided into two cases.

We first consider the case $u < 1$. We let $s = \imath \pi - \imath \epsilon e^{\imath \varphi}$, then
\begin{align*}
\mathrm{II} = e^{-\frac{d^2(u,v)}{4}} \pi^n \epsilon^{1 - n}
\int_{-\pi}^{\pi} S(n;\epsilon e^{\imath \varphi}) e^{R(u,\epsilon; \epsilon e^{\imath \varphi})} e^{\imath (1 - n )\varphi}
e^{ - \frac{\pi u}{2\epsilon}(1 - \cos{\varphi})} d\varphi,
\end{align*}
where
\begin{align*}
S(n;\xi) = \left[\frac{\xi}{\sin{\xi}} \left(1 - \frac{\xi}{\pi}\right)\right]^n.
\end{align*}
Notice that from Taylor's expansion we have
\begin{align*}
S(n;\epsilon e^{\imath \varphi}) e^{R(u,\epsilon; \epsilon e^{\imath \varphi})}
= \sum_{j = 0}^{n - 1} \eta_j(u,\epsilon) \epsilon^j e^{\imath j \varphi} + O(\epsilon^n),
\end{align*}
where $\eta_0 = 1$ and $\eta_j = O(1)$.
From \cite[p.66, p.79]{MOS66} we have
\begin{align*}
\int_{-\pi}^\pi e^{\imath (1 + j - n )\varphi}
e^{ - \frac{\pi u}{2\epsilon}(1 - \cos{\varphi})} d\varphi =
2\pi e^{ - \frac{\pi}{2\epsilon} u } I_{n - j - 1}\left(\frac{\pi}{2\epsilon} u \right),
\end{align*}
which yields
\begin{align*}
\mathrm{II} = e^{-\frac{d^2(u,v)}{4}} \pi^n \epsilon^{1 - n} \left( \sum_{j = 0}^{n - 1} 2\pi \eta_j(u,\epsilon) \epsilon^j e^{- \frac{\pi}{2\epsilon} u} I_{n - j -1}\left( \frac{\pi}{2\epsilon} u\right)
+ O(\epsilon^n)\right),
\end{align*}
where we have used $\Re \epsilon \ge |\epsilon| \cos{\eta_0^{\prime}}> 0$ in the term $O(\epsilon^n)$.

Considering separately the case $|\frac{u}{\epsilon}| < 1$ and $|\frac{u}{\epsilon}| \ge 1$ and using \eqref{I1}--\eqref{key3} as before
we have
\begin{align}\label{R2}
\frac{\eta_j(u,\epsilon) \epsilon^j I_{n - j -1}\left( \frac{\pi}{2\epsilon} u\right)}{I_{n -1}\left( \frac{\pi}{2\epsilon} u\right)} &= o(1), \\
\label{R3}
\frac{\epsilon^n}{e^{- \frac{\pi}{2\epsilon} u}I_{n -1}\left( \frac{\pi}{2\epsilon} u\right)} &= o(1),
\end{align}
which yields the result in this case. For the other case $u \ge 1$, we have
\begin{align*}
\mathrm{II} &= e^{-\frac{d^2(u,v)}{4}} \pi^n \epsilon^{1 - n}\int_{-\pi}^{\pi} e^{\imath (1 - n )\varphi}
e^{ - \frac{\pi u}{2\epsilon}(1 - \cos{\varphi})} d\varphi \\
&+ e^{-\frac{d^2(u,v)}{4}} \pi^n \epsilon^{1 - n}\int_{-\pi}^{\pi} \left( e^{R(u,\epsilon; \epsilon e^{\imath \varphi})} -1 \right)e^{\imath (1 - n )\varphi}
e^{ - \frac{\pi u}{2\epsilon}(1 - \cos{\varphi})} d\varphi \\
&+ e^{-\frac{d^2(u,v)}{4}} \pi^n \epsilon^{1 - n}\int_{-\pi}^{\pi} \left(S(n;\epsilon e^{\imath \varphi}) - 1\right) e^{R(u,\epsilon; \epsilon e^{\imath \varphi})} e^{\imath (1 - n )\varphi}
e^{ - \frac{\pi u}{2\epsilon}(1 - \cos{\varphi})} d\varphi.
\end{align*}
Note the the first integral equals to $2\pi e^{ - \frac{\pi}{2\epsilon} u } I_{n - 1}\left(\frac{\pi}{2\epsilon} u \right)$.
For the second and third integral we notice that
\begin{align*}
|R(u,\epsilon; \epsilon e^{\imath \varphi})| \le \sup_{|\xi| \le |\epsilon|}
\left|\frac{d^2}{d\xi^2} G(u;\xi)\right||\epsilon|^2 |e^{\imath \varphi} - 1|^2 \le C u |\epsilon|^2 (1 - \cos{\varphi}).
\end{align*}
Then we have
\begin{align*}
\left|\int_{-\pi}^{\pi} \left(S(n;\epsilon e^{\imath \varphi}) - 1\right) e^{R(u,\epsilon; \epsilon e^{\imath \varphi})} e^{\imath (1 - n )\varphi}
e^{ - \frac{\pi u}{2\epsilon}(1 - \cos{\varphi})} d\varphi\right| \lesssim |\epsilon| \int_{-\pi}^{\pi} e^{ - \frac{\pi u}{2|\epsilon|}
(\cos{\eta_0^{\prime}} - \frac{2C}{\pi}|\epsilon|^3)(1 - \cos{\varphi})} d\varphi.
\end{align*}
We split the second integral into the integral of $B_1 = \{|\varphi| \le \pi: C u |\epsilon| (1 - \cos{\varphi}) \le 1\}$ and
$B_2 = \{|\varphi| \le \pi: C u |\epsilon| (1 - \cos{\varphi}) > 1\}$, then
\begin{align*}
\left|\int_{B_1}\left( e^{R(u,\epsilon; \epsilon e^{\imath \varphi})} -1 \right)e^{\imath (1 - n )\varphi}
e^{ - \frac{\pi u}{2\epsilon}(1 - \cos{\varphi})} d\varphi\right| &\lesssim |\epsilon| \int_{-\pi}^{\pi} e^{ - \frac{\pi u}{2|\epsilon|}
\cos{\eta_0^{\prime}}(1 - \cos{\varphi})} d\varphi \\
\left|\int_{B_2}\left( e^{R(u,\epsilon; \epsilon e^{\imath \varphi})} -1 \right)e^{\imath (1 - n )\varphi}
e^{ - \frac{\pi u}{2\epsilon}(1 - \cos{\varphi})} d\varphi\right| &\lesssim e^{-\frac{\pi \cos{\eta_0^{\prime}}}{4C|\epsilon|^2}} \int_{-\pi}^{\pi} e^{ - \frac{\pi u}{4|\epsilon|}
(\cos{\eta_0^{\prime}} - \frac{4C}{\pi}|\epsilon|^3)(1 - \cos{\varphi})} d\varphi.
\end{align*}

Using the fact that $1 - \cos{\varphi} \ge \frac{2 \varphi^2}{\pi^2}$ for $\varphi \in [-\pi,\pi]$, we have the second and the third integral is bounded by $|\epsilon|\sqrt{\frac{|\epsilon|}{u}} = o\left(2\pi e^{ - \frac{\pi}{2\epsilon} u } I_{n - 1}\left(\frac{\pi}{2\epsilon} u \right)\right)$,
which gives the result in this case.

~ \hspace*{20pt} ~ \hfill $\Box$

\section*{Acknowledgement}
\setcounter{equation}{0}

\mbox{}\\
Ye Zhang\\
School of Mathematical Sciences  \\
Fudan University \\
220 Handan Road  \\
Shanghai 200433  \\
People's Republic of China \\
E-Mail: 17110180012@fudan.edu.cn \quad or \quad zhangye0217@126.com \mbox{}\\

\end{document}